\documentclass[10pt]{amsart}
\usepackage{geometry} 
\title{New lower bounds for the size of a non-trivial loop in the Collatz 3x+1 and generalized px+q problem}
\author{Roupam Ghosh}
\date{28 July 2009} 

\begin{document}

\maketitle
\begin{center}
\textit{Abstract}
\end{center}
In the Collatz 3x+1 problem, there are 3 possibilities: Starting from any positive number, we either reach the trivial loop (1,4,2), end up in a non-trivial loop, or go until infinity. In this paper, we shall show that if a non-trivial loop with $m$ odd numbers exists, then its minimum odd number is bounded above by a function of $m$. We shall also use that bound to calculate the least number of odd elements required for a non-trivial loop to exist. Also, the generalized bounds for the px+q problem are given.
\\
\\
\\
\\
\textbf{\Large Introduction to the Collatz problem}
\\
\\
Consider the function,
 \begin{equation*} 
f(x) = 
\begin{cases} 
\frac{x}{2} & \text{if $x$ is even}\\ 
3x+1 & \text{if $x$ is odd}
\end{cases} 
\end{equation*} 
\\
The Collatz 3x+1 conjecture states that for every positive integer $x$, there exists an integer $d (\ge 0)$ corresponding to $x$ such that 
$f^{(d)}(x) = 1$, where  $f^{(d)}(x) = f(f(f(... d $ times$ ...(f(f(f(x)))))))$
\\
For example, the iterations of $f$ on a few numbers are given below:
\begin{equation}
\begin{split}
& 1  \rightarrow 1\\
& 2 \rightarrow 1\\
& 3 \rightarrow 10 \rightarrow 5 \rightarrow 16  \rightarrow 8  \rightarrow 4  \rightarrow 2  \rightarrow 1 \\
& 4 \rightarrow 2  \rightarrow 1 \\
& 5 \rightarrow 16  \rightarrow 8  \rightarrow 4  \rightarrow 2  \rightarrow 1\\
& 6 \rightarrow  3 \rightarrow 10 \rightarrow 5 \rightarrow 16  \rightarrow 8  \rightarrow 4  \rightarrow 2 \rightarrow 1\\
& \text{and so on...}
\end{split}
\end{equation}
\\
\\
\textbf{\Large Considering only odd numbers}
\\
\\
Now, let us consider only odd integers of the Collatz sequence and modify the function $f$ to a function $T$
defined on odd integers such that $T(x) = \frac{(3x+1)}{2^k}$, where $k$ is the highest power of two in which $2^k$ divides $3x+1$.
\\
\\
For example, we have for $x = 7$:
\begin{equation}
\begin{split}
T(7) = 11,	T^{(2)}(7) = 17,	T^{(3)}(7) = 13,	T^{(4)}(7) = 5,	T^{(5)}(7) = 1,	T^{(6)}(7) = 1,	T^{(7)}(7) = 1, ...\\
k_1 = 1, k_2 = 1, k_3 = 2, k_4 = 3, k_5 = 4, k_6 = 2, k_7 = 2, ...
\end{split}
\end{equation}
\\
and so on, where $k_i$ is the highest power of two which divides $3T^{(i-1)}(x)+1$.
The iterations of $T$ on a few odd numbers are given below:
\begin{equation}
\begin{split}
& 1 \rightarrow 1\\
& 3 \rightarrow 5\rightarrow  1 \\
& 5 \rightarrow 1\\
& 7 \rightarrow 11 \rightarrow 17\rightarrow  13\rightarrow  5\rightarrow  1 \\
& 9 \rightarrow  7\rightarrow  11\rightarrow  17\rightarrow  13\rightarrow  5\rightarrow  1\\
& 11 \rightarrow 17\rightarrow  13\rightarrow  5\rightarrow  1 \\
& 13 \rightarrow  5\rightarrow  1 \\
& 15 \rightarrow 23\rightarrow  35\rightarrow  53\rightarrow  5\rightarrow  1\\
& \text{and so on...}
\end{split}
\end{equation}
\\
\\
\textbf{\Large Some inequalities}
\\
\\
Starting with any arbitrary odd integer $a_1$, let $a_r=T^{(r)}(a_1)$.
Let us suppose a non-trivial loop exists and consists of m odd numbers. Let us consider 
$2^k_1, 2^k_2...,2^k_m$ to be the highest powers of two which divide $3a_1+1, 3a_2+1, ... , 3a_m+1,$ respectively.\\
\\And let us denote $S_r = k_1 + k_2 + ... k_r$. Then we shall have
\begin{equation}
a_{m+1} = \frac{3^m}{2^{S_m}}a_{1} + \frac{c_m}{2^{S_m}} 
\end{equation}
where $c_m = 3^{m-1} + 3^{m-2} 2^{S_1} + ... + 3^2 2^{S_{m-3}} + 3^1 2^{S_{m-2}} + 2^{S_{m-1}}.$ 
\\
\bigskip
Therefore, by definition,
\begin{equation}
\begin{split}
2^{S_r} &= 2^{k_1 + k_2 + ... + k_r} \\
	     &= \frac{(3a_{1}+1)}{a_2} \frac{(3a_2+1)}{a_3} ... \frac{(3a_{r}+1)}{a_{r+1}} \\
	     &= (3+\frac{1}{a_{1}})(3+\frac{1}{a_2}) ... (3+\frac{1}{a_r}) \frac{a_{1}}{a_{r+1}} \\
	     \end{split}
\end{equation}
Then for $a_{1} = a_{m+1}$, i.e., for a loop containing $m$ odd numbers, we shall have 
\begin{equation}
2^{S_m} = (3+\frac{1}{a_{1}})(3+\frac{1}{a_2}) ... (3+\frac{1}{a_m})
\end{equation}
Let us consider $a_{min}$ to be the minimum among $a_1,a_2, ... , a_m$. Then we have \\ 
\begin{equation}
\begin{split}
& 3^m <2^{S_m} < (3+\frac{1}{a_{min}})^m \\
& m (\log_2 3) < S_m < m (\log_2 3) + m \log_2 ( 1 + \frac{1}{3a_{min}})
\end{split}
\end{equation}
\\
\\
\\
\\
\\
\\
\textbf{\Large Deriving the bound}
\\
\\
Now, since $S_m$ is a positive integer, if
$$
[m (\log_2 3) + m \log_2 ( 1 + \frac{1}{3a_{min}})]  -  [ m (\log_2 3) ]  = 0
$$
then no integer solution $S_m$ exists (where [x] is the floor function), and this will imply that no loop exists with m-odd numbers. Here, $\{x\}$ denotes $x - [x]$, the fractional part of $x$. Now, if the above condition is true then we have
\begin{equation*}
\begin{split}
& m (\log_2 3) + m \log_2 ( 1 + \frac{1}{3a_{min}}) - \{m (\log_2 3) + m \log_2 ( 1 + \frac{1}{3a_{min}})\} = m (\log_2 3) - \{ m (\log_2 3) \} \\ 
& \text{or } m \log_2 ( 1 + \frac{1}{3a_{min}}) = \{m (\log_2 3) + m \log_2 ( 1 + \frac{1}{3a_{min}})\} - \{ m (\log_2 3). \}
\end{split}
\end{equation*}
So if there is an integer solution for $S_m$, we must have
$$
m \log_2 ( 1 + \frac{1}{3a_{min}}) > \{m (\log_2 3) + m \log_2 ( 1 + \frac{1}{3a_{min}})\} - \{ m (\log_2 3) \}.
$$
This leaves us with only two possibilities, $m \log_2 ( 1 + \frac{1}{3a_{min}}) > 1$ or $m \log_2 ( 1 + \frac{1}{3a_{min}}) < 1$.
\\
\\
\\
{\textbf{\large Bound for the case $m \log_2 ( 1 + \frac{1}{3a_{min}}) > 1$:}}
\\
\\
$$
m \log_2 ( 1 + \frac{1}{3a_{min}})  > 1
$$
\\
Rearranging the terms of $m \log_2 ( 1 + \frac{1}{3a_{min}})  > 1$, we get
$$
a_{min} < \frac{1}{3 (2^\frac{1}{m}-1)}
$$
We shall define $\alpha(m)$ to be $\frac{1}{3 (2^\frac{1}{m}-1)}$. 
\\
\\
\\
{\textbf{\large Bound for the case $m \log_2 ( 1 + \frac{1}{3a_{min}})  < 1$:}}
\\
\\
From above, if a solution for $S_m$ exists then we have 
$$
m \log_2 ( 1 + \frac{1}{3a_{min}}) > \{m (\log_2 3) + m \log_2 ( 1 + \frac{1}{3a_{min}})\} - \{ m (\log_2 3) \}
$$
\\
But, since  $m \log_2 ( 1 + \frac{1}{3a_{min}})  < 1$,  $\{ m \log_2 ( 1 + \frac{1}{3a_{min}}) \} = m \log_2 ( 1 + \frac{1}{3a_{min}})$.
\\
Hence,
$$
\{ m \log_2 ( 1 + \frac{1}{3a_{min}}) \} > \{m (\log_2 3) + m \log_2 ( 1 + \frac{1}{3a_{min}})\} - \{ m (\log_2 3)\}
$$
or
$$
 \{ m (\log_2 3) \} + \{ m \log_2 ( 1 + \frac{1}{3a_{min}}) \} > \{m (\log_2 3) + m \log_2 ( 1 + \frac{1}{3a_{min}})\}$$
This is only possible if $$\{ m (\log_2 3) \} + \{ m \log_2 ( 1 + \frac{1}{3a_{min}}) \} > 1.$$
Therefore, we have
$$
m \log_2 ( 1 + \frac{1}{3a_{min}}) > 1 - \{ m (\log_2 3) \}.
$$
Rearranging the terms, we get
$$
a_{min} < \frac{1}{3 (2^{\frac{1- \{ m \log_2 3 \}}{m}}-1)}.
$$
We shall define $\beta(m)$ to be $\frac{1}{3 (2^{\frac{1- \{ m \log_2 3 \}}{m}}-1)}$.
\\
\\
\textbf{\Large Some new definitions:}
\\
\\
We shall divide the possible non-trivial loops of a Collatz sequence into two classes: We call loops satisfying the equation $$ m \log_2 ( 1 + \frac{1}{3a_{min}})  >  1 $$ $\alpha$-loops and loops satisfying the equation $$m \log_2 ( 1 + \frac{1}{3a_{min}})  <  1$$ $\beta$-loops. 
\\
Note that for $\alpha$-loops we have $ a_{min} < \alpha(m) $ and for $\beta$-loops we have $a_{min} < \beta(m)$.\\ 
\\
\\
\textbf{\Large Computing the least number of odd numbers a loop must have}
\\
\\
The Collatz algorithm has been tested and found to always reach $1$ for all numbers $\le 19 \times 2^{58}$ (Oliveira e Silva 2008).
\\
For $\alpha$-loops, we have
$$a_{min} < \frac{1}{3 (2^{\frac{1}{m}}-1)} $$
i.e.,
$$ m > \frac{\log(2)}{\log(1+\frac{1}{3 a_{min}})}
$$
Putting $a_{min} = 19 \times 2^{58} - 1$, we get
$m > 11387806137299329586$.
\\
\\
\textbf{Hence, if an $\alpha$-loop exists we must have at least $11,387,806,137,299,329,586$ odd numbers in that loop.}
\\
\\
For $\beta$-loops we have 
$$
\frac{m}{1- \{ m \log_2 3 \}} > \frac{\log(2)}{\log(1+\frac{1}{3a_{{min}}})}
$$
\\
\\
\\
\\
\\
\\
\\
\\
\\Given below is PARI/GP code I used for calculating the value of $m$ for $\beta$-loops.
\\
\begin{verbatim}
lhs(x) = ( x / (1 - frac( x*log(3.0)/log(2.0) ) ) );
rhs(x) = ( log(2)/log(1 + 1/(3*x) ) );
{
       a = rhs( 19 * ( 2^58) - 1 );
       m = 2;
       while(1,
                if( lhs(m) > a,
                         print(m," is our desired result ");
                         quit;
                         ,
                         print(m, " has been checked...");
                );		
       m++;
       );
}

\end{verbatim}
The value for $\beta$-loop is $6,586,818,670$.\\
\textbf{That is, at least $6,586,818,670$ odd numbers are required to form a $\beta$-loop, which is much
larger than the current largest known lower bound for the length of a nontrivial cycle (Sinisalo, 2003).}
\\
\\
\textbf{\Large Bounds for the generalized px+q problem}
\\
\\
Consider the function,
 \begin{equation*} 
f(x) = 
\begin{cases} 
\frac{x}{2} & \text{if $x$ is even}\\ 
px+q & \text{if } x \text{ is odd. Here } p \text{ and } q \text{ are both positive odd numbers }
\end{cases} 
\end{equation*} 
Then we shall have $$\alpha(m) = \frac{q}{p(2^\frac{1}{m}-1)}$$ $$\beta(m) = \frac{q}{p (2^{\frac{1- \{ m \log_2 p \}}{m}}-1)}$$
\\
We denote this fact by introducing this new notation,
$$\alpha_m(p,q) = \frac{q}{p(2^\frac{1}{m}-1)}$$
$$\beta_m(p,q) =  \frac{q}{p (2^{\frac{1- \{ m \log_2 p \}}{m}}-1)}$$
For example, in our case of $3x+1$, we have
$$\alpha_m(3,1) = \frac{1}{3(2^\frac{1}{m}-1)}$$
\\
\\
\\
\\
{\Large \textbf{Acknowledgements}}
\\
I thank Craig Alan Feinstein for giving me feedback on various important steps, providing me useful insights, and also amongst other things, for being a good friend. Also, my thanks goes to Eric Farin for helping me by running the program on his computer. I thank my mom, sister, and dad for making my life wonderful. 
\newline
\newline
{\Large \textbf{References}}
\begin{enumerate}
\item Weisstein, Eric W. "Collatz Problem." From MathWorld--A Wolfram Web Resource. 
 http://mathworld.wolfram.com/CollatzProblem.html 
\item Tom‡s Oliveira e Silva, http://www.ieeta.pt/\textasciitilde tos/3x+1.html
\item Shalom Eliahou (1993), "The 3x+1 Problem, New Lower Bounds on a Nontrivial Cycle Lenghts", Discrete Math 118
\item Busido Chisale (1994), "Cycles in Collatz Sequences" Publ. Math. Debrecen 45
\item Lorentz Halbeisen and Norbert Hungerbuhler (1997), "Optimal Bounds for the Length of Rationnal Collatz Cycles", Acta Arithmetica 78
\item Matti Sinisalo (2003), "On the Minimal Cycle Length of the Collatz Sequences", Univ of Oulu, Finland
\item The PARI~Group -  PARI/GP version 2.3.4, http://pari.math.u-bordeaux.fr/
\end{enumerate}
\end{document}